\documentclass[12pt]{article}
\usepackage[dvips]{graphicx}
\usepackage[T2A]{fontenc}
\usepackage[cp1251]{inputenc}
\usepackage[english]{babel}
\usepackage{babel, amsmath, amsthm, amssymb}

\begin{document}
\large

\begin{center}
{\Large A FULL PROOF OF UNIVERSAL INEQUALITIES\\ FOR THE DISTRIBUTION FUNCTION OF THE BINOMIAL LAW}\\[0.5cm]

{\large Alexander A.\,Serov, Andrew M.Zubkov\\
Steklov Mathematical Institute of RAS, Moscow, Russia}
\end{center}

\begin{quote}
{\bf Abstract:} {We present a new form and a short full proof of explicit two-sided estimates for the distribution function $F_{n,p}(x)$ of the binomial law from the paper published by D.\,Alfers and H.\,Dinges in 1984. These inequalities are universal (valid for all binomial distribution and all values of argument) and exact (namely, the upper bound for $F_{n,p}(k)$ is the lower bound for \mbox{$F_{n,p}(k+1)$).} By means of such estimates it is possible to bound any quantile of the binomial law by 2 subsequent integers.}

\end{quote}

\bigskip

Let $X_{n,p}$ be a random variable having the binomial distribution with parameters $(n,p)$:
$$
\mathbf{P}\{X_{n,p}\leqslant k\}= \sum\limits_{0\leqslant m\leqslant k} C_n^mp^m(1-p)^{n-m}.
$$
The computation of binomial sums for large $n$ being very tedious, the values of binomial distribution function are approximated usually by means of the Moivre--Laplace theorem:
\begin{equation*}
\lim\limits_{n\to\infty} \mathbf{P}\left\{X_{n,p}\leqslant np+x\sqrt{np(1-p)}\right\}=\Phi(x)\int\limits_{-\infty}^x \varphi(u)\,du,\quad
\varphi(u)=\frac1{\sqrt{2\pi}}\,e^{-u^2/2}.
\end{equation*}

For example, S.N.Bernstein \cite{ZubSer:B} proved that
$$
\sum\nolimits_{k=m_0}^{m_1-1} C_n^kp^k(1-p)^{n-k} = \frac1{\sqrt{\pi}}\, \int\nolimits_{z_0}^{z_1} e^{-u^2}du,
$$
where $m_0,\,m_1$ are integers, $z_k$ is the root of equation $z_k\sqrt{2np(1-p)}+\frac{1-2p}3\,z_k^2=m_k-np+\alpha_k$ for some $\alpha_k\in(-\tfrac32,\tfrac12)$ and $np(1-p)\geqslant 62,5,\,0\leqslant z_0<z_1\leqslant \sqrt{2np(1-p)}$; this formula is valid if $\max\{|m_0-np|,|m_1-np|\}=O\left(\sqrt{np(1-p)}\right)$. W.Feller \cite{ZubSer:F} find a modification of this formula by means of another choice of $z_0,\,z_1$. There are a lot of other results on the binomial law, see, e.\,g. \cite{ZubSer:P}. Relative errors of approximations for the tails of binomial distribution function are large due to their superexponential decreasing.

Here we present a new form and give a short full proof of some results due to D.\,Alfers and H.\,Dinges \cite{ZubSer:AD}; these authors have used some hints from \cite{PeizPratt}, \cite{Pratt}. The character of these results is analogous to the Bernstein and Feller theorems, but the formulas are explicit and (from the practical viewpoint) constitute almost final solution of the large deviation problem for the binomial law. An article \cite{ZubSer:AD} remains almost unknown (maybe because its presentation is hard to read and the proofs are too long and contain nontrivial gaps). Our proof is based on the ideas from \cite{ZubSer:AD}.

\medskip

\textbf{Theorem.} \textit{Let $H(x,p)=x\ln\tfrac{x}p+(1-x)\ln\tfrac{1-x}{1-p}$, ${\rm sgn}(x)=\tfrac{x}{|x|}$ for \mbox{$x\neq 0$} and ${\rm sgn}(0)=0$, let increasing sequences $\{C_{n,p}(k)\}_{k=0}^n$ are defined as follows $C_{n,p}(0)=(1-p)^n,\;C_{n,p}(n)=1-p^n,$
$$
C_{n,p}(k)=\Phi\left({\rm sgn}\left(\tfrac{k}n-p\right) \sqrt{2nH\left(\tfrac{k}n,p\right)}\right),\;1\leqslant k< n.
$$
Then for every $k=0,1,\ldots,n-1$ and for every $p\in(0,1)$ 
\begin{equation}
C_{n,p}(k)\leqslant\mathbf{P}\{X_{n,p}\leqslant k\}\leqslant C_{n,p}(k+1),\label{ZubSer:ADT}
\end{equation}
and equalities may happen for $k=0$ or $k=n-1$ only.}

\medskip

To demonstrate the accuracy of inequalities (\ref{ZubSer:ADT}) we may note that \mbox{$C_{n,p}(k)+C_{n,1-p}(n-k)=1$;} then from
\begin{gather*}
C_{n,p}(k)<\mathbf{P}\{X_{n,p}\leqslant k\}, \quad
C_{n,1-p}(n-k)<\mathbf{P}\{X_{n,1-p}\leqslant
n-k\}=\mathbf{P}\{X_{n,p}\geqslant k\}
\end{gather*}
it follows that
$$
1=C_{n,p}(k)+C_{n,1-p}(n-k)<\mathbf{P}\{X_{n,p}\leqslant k\}
+\mathbf{P}\{X_{n,p}\geqslant k\}=1+\mathbf{P}\{X_{n,p}=k\}.
$$
In the last inequality the difference between right and left sides is equal to the local probability of the binomial law. So, $\mathbf{P}\{X_{n,p}\leqslant k\}-C_{n,p}(k)<\mathbf{P}\{X_{n,p}= k\}$.

The ratio of upper and lower bounds for $\mathbf{P}\{X_{n,p}\leqslant k\}$ in (\ref{ZubSer:ADT}) may be large if $k$ is significantly less than $np$, but for such $k$ the ratios  $\mathbf{P}\{X_{n,p}\leqslant k+1\}/\mathbf{P}\{X_{n,p}\leqslant k\}$ are large also.

In somewhat another form we use the results of \cite{ZubSer:AD} in \cite{ZubSer:ZS} to estimate the partial sums of binomial coefficients.

\medskip

\textsl{Proof.} Lower bound in (\ref{ZubSer:ADT}) for $k=0$ and upper bound for $k=n-1$ are exact equalities; so these cases will not be considered further.

We have \mbox{$\mathbf{P}\{X_{n,1}\leqslant k\}=0,k<n$,} and for every integer $k\in\{0,1,\ldots,n-1\}$
\begin{gather*}
\mathbf{P}\{X_{n,p}\leqslant k\}= \sum\limits_{m=0}^k C_n^m
p^m(1-p)^{n-m}= -\int_p^1 \frac{d}{dz}\,\sum\limits_{m=0}^k C_n^m
z^m(1-z)^{n-m}\,dz= \\
=-n\int_p^1 \sum\limits_{m=0}^k
(C_{n-1}^{m-1} z^{m-1}(1-z)^{n-m}-C_{n-1}^m z^m(1-z)^{n-m-1})\,dz=\\
=n C_{n-1}^k\int_p^1
z^k(1-z)^{n-k-1}dz=(k+1)C_n^{k+1} \int_p^1\frac{z^{k+1}(1-z)^{n-k-1}}{z}\,dz.
\end{gather*}

Applying the Stirling formula $n!=\sqrt{2\pi
n}\left(\frac{n}e\right)e^{S_n}$ we find that
\begin{gather*}
(k+1)C_n^{k+1}=\exp\{S_n^{k+1}\}\sqrt{\frac{(k+1)n}{2\pi(n-k-1)}}\,\frac{n^n}{(k+1)^{k+1}(n-k-1)^{n-k-1}}\,,
\end{gather*}
where $S_n^{k+1}=S_n-S_{k+1}-S_{n-k-1}$. 

The main step is the proof of the upper bound in (\ref{ZubSer:ADT}) for $0\leqslant k\leqslant n-2$. Denoting $\alpha=\frac{k+1}n\in(0,\,1)$ we have
\begin{gather*}
\frac{n^n}{(k+1)^{k+1}(n-k-1)^{n-k-1}}\,z^{k+1}(1-z)^{n-k-1}= \left(\frac{nz}{k+1}\right)^{k+1}\left(\frac{n(1-z)}{n-k-1}\right)^{n-k-1}=\\
=\left(\frac{z}\alpha\right)^{n\alpha}\left(\frac{1-z}{1-\alpha}\right)^{n(1-\alpha)}
=\exp\left\{-n\left(\alpha\ln\frac{\alpha}z+(1-\alpha)\ln\frac{1-\alpha}{1-z}\right)\right\}.
\end{gather*}
So, if $B(z)\stackrel{\text{def}}{=}\alpha\ln\frac{\alpha}{z}+
(1-\alpha)\ln\frac{1-\alpha}{1-z}$ then
\begin{gather*}
\mathbf{P}\{X_{n,p}\leqslant k\}=e^{S_n^{k+1}}\sqrt{\frac{n\alpha}{2\pi(1-\alpha)}}
\int_p^1 e^{-nB(z)}\,\frac{dz}z\,.\label{X>=k0}
\end{gather*}
For each $\alpha\in(0,\,1)$ the function $B(z)$ decreases monotonically from $+\infty$ to 0 on $(0,\,\alpha]$ and increases monotonically from 0 to $+\infty$ on $[\alpha,\,1)$, indeed:
\begin{gather}
B(\alpha)=0\quad \text{and}\quad B'(z)=-\frac{\alpha}z+\frac{1-\alpha}{1-z}=\frac{z-\alpha}{z(1-z)}\,. \label{B'}
\end{gather}
It follows that the equation $B(z)=\frac12\,a^2(z)$ has solution $a(z)=(\alpha-z)\sqrt{\frac{2B(z)}{(\alpha-z)^2}}$ such that $a(z)$ decreases monotonically on $(0,\,1)$ from $+\infty$ to $-\infty$. So,
\begin{gather}
\mathbf{P}\{X_{n,p}\leqslant k\}=\nonumber\\
=e^{S_n^{k+1}}\sqrt{\frac{n\alpha}{2\pi(1-\alpha)}}
\int_p^1 e^{-n\frac{a^2(z)}2}\,\frac{dz}z
=e^{S_n^{k+1}}\sqrt{\frac{n\alpha}{1-\alpha}} \int_p^1 \varphi\left(a(z)\sqrt{n}\right)\frac{dz}z\,.\label{X<=k}
\end{gather}

Further, let us consider the difference
$$
\delta(p)=\mathbf{P}\{X_{n,p}\leqslant k\}-\Phi\left(a(p)\sqrt{n}\right),\quad p\in[0,\,1],
$$
where $\Phi(\cdot)$ is the standard normal distribution function. If $0\leqslant k\leqslant n-2$ then $\delta(0)=\delta(1)=0$. Now to prove that $\delta(p)<0$ for all $p\in(0,\,1)$ it is sufficient to show that: a) the function $\delta(p)$ is differentiable with respect to $p$, b) the equation $\delta'(p)=0$ has unique root $p_0$ on $(0,\,1)$ and c) $\delta'(p)<0$ for $p\in(0,\,p_0)$ and $\delta'(p)>0$ for $p\in(p_0,\,1)$.

In view of (\ref{X<=k})
$$
\frac{d}{dp}\,\mathbf{P}\{X_{n,p}\leqslant k\} = -e^{S_n^{k+1}}\sqrt{\frac{n\alpha}{1-\alpha}}\,\frac1p\, \varphi\left(a(p)\sqrt{n}\right).
$$
Further, it follows from $B(p)=\frac12\,a^2(p)$ and (\ref{B'}) that $B'(p)=\frac{p-\alpha}{p(1-p)}=a(p)a'(p)$, i.\,e. $a'(p)=\frac{p-\alpha}{p(1-p)a(p)}$\,, and so
\begin{gather*}
\frac{d}{dp}\,\Phi\left(a(p)\sqrt{n}\right) = \varphi\left(a(p)\sqrt{n}\right)a'(p)\sqrt{n}= \varphi\left(a(p)\sqrt{n}\right)\frac{\sqrt{n}(p-\alpha)}{p(1-p)a(p)}=\\
=-\varphi\left(a(p)\sqrt{n}\right)\frac{\sqrt{n}}{p(1-p)}
\sqrt{\frac{(p-\alpha)^2}{2B(p)}}\,.
\end{gather*}
We have
\begin{gather}
\delta'(p)=-\frac1p\, \varphi\left(a(p)\sqrt{n}\right)\sqrt{n}\left(e^{S_n^{k+1}}\sqrt{\frac{\alpha}{1-\alpha}} - \sqrt{\frac{(p-\alpha)^2}{2(1-p)^2B(p)}}\right).
\end{gather}
The multiplier before parenthesis is negative, the first term in parenthesis doesn't depend on $p$. From the formula
$B(p)=\alpha\ln\left(\frac{\alpha}p\right)+ (1-\alpha)\ln\left(\frac{1-\alpha}{1-p}\right)$ it follows that
$$
\lim\limits_{p\downarrow 0}\sqrt{\frac{(p-\alpha)^2}{(1-p)^2B(p)}}=0,\quad \lim\limits_{p\uparrow 1}\sqrt{\frac{(p-\alpha)^2}{(1-p)^2B(p)}}=\infty.
$$
Let us show that the function $\frac{(p-\alpha)^2}{2(1-p)^2B(p)}$ is monotonically increasing. Its derivative equals to
\begin{gather*}
\frac{\partial}{\partial p}\frac{(p-\alpha)^2}{(1-p)^2B(p)}= \frac1{B(p)}\,\frac{\partial}{\partial p}\frac{(p-\alpha)^2}{(1-p)^2}+\frac{(p-\alpha)^2}{(1-p)^2}\,\frac{\partial}{\partial p}\frac1{B(p)}= \\
=2\frac{(p-\alpha)(1-\alpha)}{(1-p)^3B(p)}- \frac{(p-\alpha)^2}{(1-p)^2}\,\frac{p-\alpha}{p(1-p)B^2(p)}=\\
=\frac{2(1-\alpha)(p-\alpha)}{(1-p)^3B^2(p)}\,\left(B(p)-\frac{(p-\alpha)^2}{2p(1-\alpha)}\right).
\end{gather*}
The first multiplier in the right hand side changes its sign from $-$ to $+$ at $p=\alpha$, the difference in the parenthesis equals to 0 for $p=\alpha$, and
$$
\frac{\partial}{\partial p}\left(B(p)-\frac{(p-\alpha)^2}{2p(1-\alpha)}\right)= \frac{p-\alpha}{p(1-p)}-\frac{p^2-\alpha^2}{2p^2(1-\alpha)}=\frac{(p+1)(p-\alpha)^2}{2p^2(1-p)(1-\alpha)} > 0
$$
for $p\in(0,1)\backslash\{\alpha\}$, i.\,e. the difference $B(p)-\frac{(p-\alpha)^2}{2p(1-\alpha)}$ changes its sign from $-$ to $+$ at $p=\alpha$ also. It means that $\frac{\partial}{\partial p}\frac{(p-\alpha)^2}{(1-p)^2B(p)}> 0$ for all $p\in(0,1)\backslash\{\alpha\}$.

Consequently, the difference $e^{S_n^k}\sqrt{\frac{\alpha}{1-\alpha}} - \sqrt{\frac{(p-\alpha)^2}{2(1-p)^2B(p)}}$ on $[0,\,1]$ decreases monotonically from $e^{S_n^k}\sqrt{\frac{\alpha}{1-\alpha}}$ to $-\infty$ and is equal to 0 for some $p=p_0$. We see that $\delta(p)< 0$ for all $p\in(0,\,1)$, i.\,e. for all $n\geqslant 1,\, 0\leqslant k\leqslant n-2,\,p\in(0,\,1)$
\begin{gather*}
\mathbf{P}\{X_{n,p}\leqslant k\}< \Phi\left(a(p)\sqrt{n}\right)= \label{X<kfin>}\\
=\Phi\left({\rm sgn}\left(\tfrac{k+1}n-p\right)\sqrt{2n \left(\tfrac{k+1}n\,\ln\tfrac{k+1}{np}+\left(1-\tfrac{k+1}n\right)\ln\tfrac{n-k-1}{n(1-p)}\right)}\right)=\nonumber\\
=\Phi\left({\rm sgn}\left(\tfrac{k+1}n-p\right)\sqrt{2n H\left(\tfrac{k+1}n\,,p\right)}\right)=C_{n,p}(k+1).\nonumber
\end{gather*}

Now we prove the lower bound in (\ref{ZubSer:ADT}) for $k\in\{1,\ldots,n-1\}$. As
$$
\mathbf{P}\{X_{n,p}\leqslant k\}=\mathbf{P}\{X_{n,1-p}\geqslant n-k\}=1-\mathbf{P}\{X_{n,1-p}\leqslant n-k-1\}
$$
and $n-k-1\in\{0,\ldots,n-2\}$ we have (using upper bound just proved)
\begin{gather*}
\mathbf{P}\{X_{n,p}\leqslant k\}> 1-\Phi\left({\rm sgn}\left(\tfrac{n-k}n-(1-p)\right)\sqrt{2n H\left(\tfrac{n-k}n\,,1-p\right)}\right)=\nonumber\\
=\Phi\left({\rm sgn}\left(\tfrac{k}n-p\right)\sqrt{2n H\left(\tfrac{k}n\,,p\right)}\right)=C_{n,p}(k).\label{X<kfin>}
\end{gather*}
The theorem is proved.

\medskip

\textsl{Remark.} Inequalities (\ref{ZubSer:ADT}) may be sharpened by means of nontrivial upper bounds of the function
\begin{gather*}
\delta(p)=e^{S_n^{k+1}}\sqrt{\frac{n\alpha}{2\pi(1-\alpha)}} \int_p^1
e^{-n\frac{a^2(z)}2}\,\frac{dz}z-\Phi\left(a(p)\sqrt{n}\right)=\\
=-\int\nolimits_0^p \frac1z\, \varphi\left(a(z)\sqrt{n}\right)\sqrt{n}\left(e^{S_n^{k+1}}\sqrt{\frac{\alpha}{1-\alpha}} - \sqrt{\frac{(z-\alpha)^2}{2(1-z)^2B(z)}}\right)dz,\quad p\in(0,\,1).
\end{gather*}
For concrete values of parameters this integral may be estimated numerically.

To obtain the analytic estimates we may note that the integrand is the product of continuous positive function
$f(z)=\frac1z\, \varphi\left(a(z)\sqrt{n}\right)\sqrt{n}$ and monotonically decreasing function
\mbox{$g(z)=e^{S_n^{k+1}}\sqrt{\frac{\alpha}{1-\alpha}} -
\sqrt{\frac{(z-\alpha)^2}{2(1-z)^2B(z)}}$, and that
$\int\nolimits_0^1 f(z)g(z)dz=0$.} So, if $p_0\in(0,1)$ is such that
\mbox{$g(p_0)=0$,} then the value of the integral may be bounded from below: for $0<p\leqslant p_0$ we have
$$
\int\nolimits_0^p f(z)g(z)dz\geqslant \sup\limits_{0<u\leqslant p}g(u)\int\nolimits_0^u f(z)dz,
$$
and for $p_0\leqslant p<1$
$$
\int\nolimits_0^p f(z)g(z)dz=-\int\nolimits_p^1 f(z)g(z)dz\geqslant \sup\limits_{p\leqslant u<1}|g(u)|\int\nolimits_u^1 f(z)dz.
$$
These estimates along with (\ref{X<=k}) and (\ref{ZubSer:ADT}) give for $p\geqslant p_0$
\begin{gather}
\delta(p)=-\int\nolimits_p^1 \delta'(y)\,dy \leqslant\nonumber\\
=-\sup\limits_{z\in(p,1)}\left|e^{S_n^{k+1}}\sqrt{\frac{\alpha}{1-\alpha}} - \sqrt{\frac{(z-\alpha)^2}{2(1-z)^2B(z)}}\right|\int\nolimits_z^1 \frac1y\, \varphi\left(a(y)\sqrt{n}\right)\sqrt{n}\,dy =\nonumber\\
=-\sup\limits_{z\in(p,1)}\left|1-e^{-S_n^{k+1}} \sqrt{\frac{(1-\alpha)(z-\alpha)^2}{2\alpha(1-z)^2B(z)}}\right|\mathbf{P}\{X_{n,z}\leqslant k\}\leqslant\nonumber\\
\leqslant -\sup\limits_{z\in(p,1)}\left|1-e^{-S_n^{k+1}} \sqrt{\frac{(1-\alpha)(z-\alpha)^2}{2\alpha(1-z)^2B(z)}}\right|C_{n,z}(k).\label{delta+}
\end{gather}
Analogously, for $p<p_0$
\begin{gather}
\delta(p)=\int\nolimits_0^p \delta'(y)\,dy \leqslant\nonumber\\
=-\sup\limits_{z\in(0,p)}\left|e^{S_n^{k+1}}\sqrt{\frac{\alpha}{1-\alpha}} - \sqrt{\frac{(z-\alpha)^2}{2(1-z)^2B(z)}}\right|\int\nolimits_0^z \frac1y\, \varphi\left(a(y)\sqrt{n}\right)\sqrt{n}\,dy =\nonumber\\
=-\sup\limits_{z\in(0,p)}\left|1-e^{-S_n^{k+1}} \sqrt{\frac{(1-\alpha)(z-\alpha)^2}{2\alpha(1-z)^2B(z)}}\right|\mathbf{P}\{X_{n,z}\geqslant k+1\}\leqslant\nonumber\\
\leqslant -\sup\limits_{z\in(0,p)}\left|1-e^{-S_n^{k+1}} \sqrt{\frac{(1-\alpha)(z-\alpha)^2}{2\alpha(1-z)^2B(z)}}\right|(1-C_{n,z}(k+1)).\label{delta-}
\end{gather}


\begin{thebibliography}{4}

\bibitem{ZubSer:AD} {\it Alfers D., Dinges H.}, A normal approximation for
Beta and Gamma tail probabilities. --- Z. Wahrscheinlichkeitstheor. verw. Geb., 1984, v.~65, p.~399--420.

\bibitem{ZubSer:B} {\it Bernstein S.\,N.}, Retour au probl\`{e}me de l'\'{e}valuation de l'approximation de la formule limite de Laplace (in Russian). --- Izvestia Akademii Nauk SSSR, ser. matem., 1943, v.~7, №~1, p.~3--16.

\bibitem{ZubSer:F} {\it Feller W.}, On the normal approximation to the binomial distribution. ---
Ann. Math. Statist., 1945, v.~16, №~4, p.~319--329.

\bibitem{ZubSer:P} {\it Petrov V.\,V.}, Limit theorems of probability theory. Sequences of independent random variables. ---  Clarendon Press, 1995.

\bibitem{PeizPratt} {\it Peizer D.\,B., Pratt J.\,W.} A normal approximation for binomial, $F$, beta and other common related tail probabilities (part I). --- J. Amer. Statist. Ass., 1968, no. 63, p. 1416-1456.

\bibitem{Pratt} {\it Pratt J.\,W.} A normal approximation for binomial, $F$, beta and other common related tail probabilities (part II). --- J. Amer. Statist. Ass., 1968, no. 63, p. 1457-1483.

\bibitem{ZubSer:ZS} {\it Zubkov A.\,M, Serov A.\,A.}, Bounds for the number of Boolean functions 
admitting affine approximations of a given accuracy. --- Discrete Math. Appl., 2010, 20, №~5-6, p.~467--486.

\end{thebibliography}
\end{document}